\documentclass[righttag,12pt]{amsart}
\usepackage{amssymb}
\usepackage{euscript}

 \let\cal\mathcal

\newtheorem{prop}{Proposition}[section]
\newtheorem{theorem}[prop]{Theorem}
\newtheorem{cor}[prop]{Corollary}
\newtheorem{lem}[prop]{Lemma}
 
\numberwithin{equation}{section}

\theoremstyle{definition}

\newtheorem{ack}{Acknowledgments}

\theoremstyle{remark}
\newtheorem{rem}[prop]{Remark}
\newcommand{\bbN}{{\mathbb{N}}}
\def \eps{\varepsilon}

\def\Ndb{\mathbb N}

\def\Rdb{\mathbb R}

\newcommand{\al}{\alpha}

\newcommand{\de}{\delta}
\newcommand{\e}{\varepsilon}

\newcommand{\ro}{\varrho}

\newcommand{\f}{\varphi}

\newcommand{\sgn}{\operatorname{sgn}}

 \newcommand{\disp}{\displaystyle}

\newcommand{\lb}{\label}

\newcommand{\lra}{\longrightarrow}
 
\begin{document}

\title{On the extension of H\"{o}lder maps with values in
spaces of continuous functions}

\author{Gilles Lancien}
\address{D\'{e}partement de Math\'{e}matiques\\ Universit\'{e} de Franche-Comt\'{e}\\
16 Route de Gray, 25030 Besan\c{c}on, France}

\email{glancien@math.univ-fcomte.fr}

\author{Beata Randrianantoanina$^*$}\thanks{$^*$Partially
funded by a CFR Grant from Miami University}

\address{Department of Mathematics and Statistics \\ Miami University
\\Oxford, OH 45056, USA}

\email{randrib@muohio.edu}

\begin{abstract} We study the isometric extension problem for H\"{o}lder
maps from subsets of any Banach space into $c_0$ or into a space
of continuous functions. For a Banach space $X$, we prove that any
$\alpha$-H\"{o}lder map, with $0<\alpha\leq 1$, from a subset of $X$
into $c_0$ can be isometrically extended to $X$ if and only if $X$
is finite dimensional. For a finite dimensional normed space $X$
and for a compact metric space $K$, we   prove that the set of
$\alpha$'s for which all $\alpha$-H\"{o}lder maps from a subset of $X$
into $C(K)$ can be extended isometrically is either $(0,1]$ or
$(0,1)$ and we give examples of both occurrences. We also prove
that for any metric space $X$, the described above set of $\al$'s
does not depend on $K$, but only on finiteness of $K$.

\end{abstract}


\maketitle

\section{Introduction - Notation}

If $(X,d)$ and $(Y,\ro)$ are metric spaces, $\alpha \in (0,1]$ and
$K>0$ , we will say that a map $f:X \to Y$ is $\alpha$-H\"{o}lder with
constant $K$ (or in short $(K,\alpha)$-H\"{o}lder) if
$$\forall x,y \in X,\ \ \ro(f(x),f(y))\leq Kd(x,y)^\alpha.$$
Let us now recall and extend the notation introduced by Naor in
\cite{N}. For $C\geq 1$, ${\cal B}_C(X,Y)$ will denote the set of
all $\alpha \in (0,1]$ such that any $(K,\alpha)$-H\"{o}lder function
$f$ from a subset of $X$ into $Y$ can be extended to a
$(CK,\alpha)$-H\"{o}lder function from $X$ into $Y$. If $C=1$, such an
extension is called an isometric extension. When $C>1$, it is
called an isomorphic extension. If a $(CK,\alpha)$-H\"{o}lder
extension exists for all $C>1$, we will say that $f$ can be almost
isometrically extended. So, let us define:
$${\cal A}(X,Y)={\cal B}_1(X,Y),\ \ {\cal B}(X,Y)=\bigcup_{C\geq 1}{\cal B}_C(X,Y),
\ {\rm and}\ \widetilde{\cal A}(X,Y)=\bigcap_{C> 1}{\cal
B}_C(X,Y).$$

The study of these sets goes back to a classical result of
Kirszbraun \cite{K} asserting that if $H$ is a Hilbert space, then
$1 \in {\cal A}(H,H)$. This was extended by Gr\"{u}nbaum and
Zarantonello \cite {GZ} who showed that ${\cal A}(H,H)=(0,1]$.
Then the complete description of ${\cal A}(L^p,L^q)$ for
$1<p,q<\infty$ relies on works by Minty \cite{M} and Hayden, Wells
and Williams \cite{HWW} (see also the book of Wells and Williams
\cite{WW} for a very nice exposition of the subject). More
recently, K. Ball \cite{B} introduced a very important notion of
non linear type or cotype and used it to prove a general extension
theorem for Lipschitz maps. Building on this work, Naor \cite{N}
described completely the sets ${\cal B}(L^p,L^q)$ for
$1<p,q<\infty$.

In this paper, we concentrate on the study of ${\cal A}(X,Y)$ and
$\widetilde{\cal A}(X,Y)$, when $X$ is a Banach space and $Y$ is a
space of converging sequences or, more generally, a space of
continuous functions on a compact metric space. This can be viewed
as an attempt to obtain a non linear version of the results of
Lindenstrauss and Pelczy\'{n}ski \cite{LP} and later of Johnson
and Zippin (\cite{JZ1} and \cite{JZ2}) on the extension of linear
operators with values in $C(K)$ spaces.

\noindent So let us denote by $c$ the space of all real converging
sequences equipped with the supremum norm and by $c_0$ the
subspace of $c$ consisting of all sequences converging to 0. If
$K$ is a compact space, $C(K)$ denotes the space of all real
valued continuous functions on $K$, equipped again with the
supremum norm.

\noindent In section 2, we show that if $X$ is infinite
dimensional and $Y$ is any separable Banach space containing an
isomorphic copy of $c_0$, then $\widetilde{\cal A}(X,Y)$ is empty.
On the other hand, we prove that ${\cal A}(X,c_0)=(0,1]$, whenever
$X$ is finite dimensional.

\noindent In section 3, we show that for any finite dimensional
space $X$, $\widetilde{\cal A}(X,c)=(0,1]$ and ${\cal A}(X,c)$
contains $(0,1)$. Then the study of the isometric extension for
Lipschitz maps turns out to be a bit more surprising. Indeed, we
give an example of a 4-dimensional space $X$ such that ${\cal
A}(X,c)=(0,1)$. To our knowledge, this provides the first example
of Banach spaces $X$ and $Y$ such that ${\cal A}(X,Y)$ is not
closed in $(0,1]$  and also such that ${\cal A}(X,Y)\neq
\widetilde{\cal A}(X,Y)$. On the other hand, we show that if the
unit ball of a finite dimensional Banach space is a polytope, then
${\cal A}(X,c)=(0,1]$.

\noindent Finally, we prove in section 4, that $c$ is the only
$C(K)$ space that one needs to consider as the image space in the
study of the isometric extension problem. More precisely, we show
that for every infinite compact metric space $K$ and every metric
space $X$, ${\cal A}(X,c)={\cal A}(X,C(K))$.

\begin{ack} The research on this paper started during a sabbatical
visit of the second named author at the D\'{e}partement de
Math\'{e}matiques, Universit\'{e}
 de Franche-Comt\'{e} in Besan\c{c}on, France. She wishes to thank all members of
the Functional Analysis Group, and especially Prof. F. Lancien,
 for their hospitality during that visit.
\end{ack}

\section{Maps into $c_0$ }
It is well known that for any metric space $(X,d)$, ${\cal
A}(X,\Rdb)=(0,1]$. Indeed, if $M$ is a subset of $X$ and $f:M \to
\Rdb$ is a $(K,\alpha)$-H\"{o}lder function, then a
$(K,\alpha)$-H\"{o}lder extension $g$ of $f$ on $X$ is given for
instance by the inf-convolution formula:
$$\forall x \in X,\ \ g(x)=\inf\{f(u)+K(d(u,x))^\alpha,\ u\in
M\}.$$ It follows immediately that ${\cal
A}(X,\ell_\infty)=(0,1]$, where $\ell_\infty$ is the space of all
real bounded sequences equipped with the supremum norm . Now,
since there is a 2-Lipschitz retraction from $\ell_\infty$ onto
$c_0$ (see for instance \cite[page 14]{BL}), it is clear that for
any metric space $X$, ${\cal B}_2(X,c_0)=(0,1]$. Our first result
shows that the difference between the isometric and isomorphic
extension problems which is revealed in \cite{N} is extreme when
$c_0$ is the image space. More precisely:

\begin{theorem}\lb{infdim} Let $X$ be an infinite dimensional
normed vector space and $Y$ be a separable Banach space containing
an isomorphic copy of $c_0$. Then
$${\widetilde A}(X,Y)=\emptyset.$$
\end{theorem}

\begin{proof} By a theorem of R.C. James \cite{J}, $Y$
contains almost isometric copies of $c_0$. So, since we are
studying the almost isometric extension problem, we may as well
assume that there is a closed subspace $Z$ of $Y$ which is
isometric to $c_0$. Let $(e_n)$ be the isometric image in $Z$ of
the canonical basis of $c_0$ and $(e_n^*)$ be the Hahn-Banach
extensions to $Y$ of the corresponding coordinate functionals
(this sequence is included in the unit sphere of $Y^*$). Since $Y$
is separable, there is a subsequence $(e^*_{n_k})_{k\geq 1}$ which
is weak$^*$-converging to some $y^*$ in the unit ball of $Y^*$.

\noindent On the other hand, by a theorem of Elton and Odell
\cite{EO}, there exists $\eps>0$ and a sequence $(x_k)_{k\geq 1}$
in $X$ such that:
$$\forall k\ \ \Vert x_k\Vert=1-\eps\ \ {\rm and}\ \ \forall k\neq
l\ \ \Vert x_k-x_l\Vert\geq 1.$$

\noindent Let now $f$ be defined by $f(x_k)=(-1)^ke_{n_k}$. This
is clearly a $(1,\alpha)$-H\"{o}lder function for any $\alpha$ in
$(0,1]$. Let $\delta >0$ such that $(1+\delta)(1-\eps)^\alpha <1$
and $\eta=1-(1+\delta)(1-\eps)^\alpha>0$. Assume that $f$ can be
extended at 0 into a $(1+\delta,\alpha)$-H\"{o}lder function $g$ with
$g(0)=y$. Then, for any even $k$, $e^*_{n_k}(y)\geq \eta$ and for
any odd $k$, $e^*_{n_k}(y)\leq -\eta$. This is in contradiction
with the fact that $(e^*_{n_k})$ is weak$^*$-converging.

\end{proof}

We will now  solve the extension problem for H\"{o}lder maps from a
finite dimensional space into $c_0$. First, we need the following
elementary Lemma.
\begin{lem}\lb{cones} Let $X$ be a finite dimensional Banach space and
$\delta>0$. Then there exist $C_1,...,C_n$ subsets of $X$ such
that
$$X\setminus \{0\}=\bigcup_{i=1}^n C_i$$
and
$$\forall 1\leq i\leq n,\ \forall x,y \in C_i\ {\rm so \ that}\
\Vert x\Vert\geq \Vert y\Vert:\ \Vert x-y\Vert\leq \Vert x\Vert
-(1-\delta)\Vert y\Vert.$$
\end{lem}

\begin{proof} Since $X$ is finite dimensional, we can cover
the unit sphere of $X$ with $B_1,...,B_n$, balls of radius
$\delta/2$ and define
$$C_i=\{y\in X\setminus \{0\}: \ {y \over \Vert y\Vert}\in
B_i\}.$$ Let now $x,y \in C_i$ so that $\Vert x\Vert\geq \Vert
y\Vert.$ We have

\begin{equation*}
\begin{split}
\Vert x-y\Vert \leq \big\Vert x-x\frac{\Vert y\Vert}{\Vert
x\Vert}\big\Vert +\Vert y\Vert \,\big\Vert \frac{x}{\Vert
x\Vert}-\frac{y}{\Vert y\Vert}\big\Vert \leq \Vert x\Vert -
(1-\delta)\Vert y\Vert.
\end{split}
\end{equation*}

\end{proof}

Then our result is

\begin{theorem}\lb{fdc0} If $X$ is a finite dimensional normed vector
space, then $${\cal A}(X,c_0)=(0,1].$$
\end{theorem}

\begin{proof} Let $\alpha \in (0,1]$, $M\subset X$ and $f:M\to
c_0$ be a $(K,\alpha)$-H\"{o}lder function. We may assume that $K=1$
and that $M$ is closed. It is enough to show that for any $x_0 \in
X\setminus M$, $f$ can be extended into a $(1,\alpha)$-H\"{o}lder
function $g$ on $M\cup \{x_0\}$ and we will assume that $x_0=0$.

For $\delta=1/2$, let $C_1,..,C_n$ be given by Lemma \ref{cones}
and $I=\{i,\ 1\leq i\leq n\ {\rm and}\ C_i\cap M \neq
\emptyset\}.$ Since $X$ is finite dimensional, for each $i$ in
$I$, we can pick $x_i$ in $\overline{C_i}\cap M$ such that for any
$x\in C_i\cap M$, $\Vert x\Vert\geq \Vert x_i\Vert$. Then, by
Lemma \ref{cones}, we have that

\begin{equation*}
\forall x\in C_i\cap M,\ \Vert x-x_i\Vert\leq \|x\|-\frac12\|x_i\|
\leq \Vert x\Vert.
\end{equation*}

\noindent Let us now pick $\eps >0$ such that $\eps<{\displaystyle
\frac12} dist(0,M)^\alpha.$ Then

\begin{equation*}
\exists N\in \Ndb\ \forall n>N\ \forall i\in\{1,..,n\}:\  \vert
f(x_i)(n)\vert<\eps.
\end{equation*}

\noindent We will now choose $g(0)=(u(n))_{n\geq 1}$.

\noindent Since $\Rdb$-valued contractions can be extended into
contractions, we can pick $(\eta_n)_{n\geq 1}$ in $\ell_\infty$ so
that
$$\forall n\in \Ndb \ \forall x \in M,\ \vert f(x)(n)-\eta_n\vert \leq
\Vert x\Vert^\alpha.$$

\noindent For $n\leq N$, we set $u(n)=\eta_n$.

\noindent For $n>N$, let $\delta_n \in \{-1,1\}$ be the sign of
$\eta_n$. Now we set
$$u(n)= \de_n \min\{|\eta_n|, \max_{i\in I} |f(x_i)(n)|\}.$$
Note that since $I$ is finite and each $f(x_i)\in c_0$, we have
that $g(0)=(u(n))_{n\geq 1} \in c_0$.

\noindent Next we check that for all $x\in M$ and all $n>N$,
$\vert f(x)(n)-u(n)\vert \leq \Vert x\Vert^\alpha.$ So let $x\in
M$ and $i_0\in I$ such that $x\in C_{i_0}\cap M$. We have four
cases:

1) If $|f(x)(n)|\le |u(n)|$, then
$$|f(x)(n)-u(n)|\le 2\e\le \|x\|^\al.$$

2) If $|f(x)(n)|> |u(n)|$, $\sgn (f(x)(n))=\de_n$, and
$|u(n)|=|\eta_n|$ then, by the definition of $\eta_n$:
$$|f(x)(n)-u(n)|\le \|x\|^\al.$$

3) If $|f(x)(n)|> |u(n)|$, $\sgn (f(x)(n))=\de_n$, and
$|u(n)|=\disp{\max_{i \in I}} |f(x_i)(n)|\ge |f(x_{i_0})(n)|$,
then
\begin{equation*}
\begin{split}
|f(x)(n)-u(n)|&= |f(x)(n)|-|u(n)|\le  |f(x)(n)-f(x_{i_0})(n)|+ |f(x_{i_0})(n)|-|u(n)|\\
&\le \|x-x_{i_0}\|^\al \le \|x\|^\al.
\end{split}
\end{equation*}

4) If $|f(x)(n)|> |u(n)|$ and $\sgn (f(x)(n))\ne\de_n$, then
\begin{equation*}
|f(x)(n)-u(n)|= |f(x)(n)|+|u(n)|\le
|f(x)(n)|+|\eta_n|=|f(x)(n)-\eta_n|\le\|x\|^\al.
\end{equation*}

\end{proof}

\begin{rem}  The proof is much simpler in the case
$\alpha=1$. Indeed it is enough to set $u(n)=0$ for $n>N$. Then,
for $x\in M$, pick $i_0\in I$ such that $x\in C_{i_0}\cap M$.
Thus, for all $n>N$:

\begin{equation*}
\begin{split}
\vert f(x)(n)-u(n)\vert=\vert f(x)(n)\vert&\leq \vert
f(x)(n)-f(x_{i_0})(n)\vert+\e\leq \Vert x-x_{i_0}\Vert+\e \\ &\leq
\Vert x\Vert-\frac12 \Vert x_{i_0}\Vert +\e \leq \Vert x\Vert.
\end{split}
\end{equation*}
\end{rem}

\section{Maps into $c$}

We now consider the isometric and almost isometric extension
problems for H\"{o}lder maps from a normed vector space into $c$. If
$X$ is infinite dimensional, this question is settled by Theorem
\ref{infdim}. Therefore, throughout this section, $X$ will denote
a finite dimensional normed vector space. The study of the almost
isometric extensions is then rather simple. For this purpose, we
recall that, for $\lambda >1$, a Banach space $Y$ is said to be a
${\cal L}^\infty_\lambda$ space if every finite dimensional
subspace of $Y$ is contained in a finite dimensional subspace $F$
of $Y$ which is $\lambda$-isomorphic to $\ell_\infty^{{\rm dim}\,
F}$ (namely, there is an isomorphism $T$ from $F$ onto
$\ell_\infty^{{\rm dim}\, F}$ such that $\Vert T\Vert\,\Vert
T^{-1}\Vert\leq \lambda$).

\begin{prop}\lb{almost} Let $X$ be a finite dimensional normed vector space
and $Y$ be a Banach space which is a ${\cal L}^\infty_\lambda$
space for any $\lambda>1$. Then $${\widetilde A}(X,Y)=(0,1].$$ In
particular, for every compact space $K$, $${\widetilde
A}(X,C(K))=(0,1].$$
\end{prop}

\begin{proof} Let $M$ be a closed subset of $X$ and $f:M\to Y$ be a
$(1,\alpha)$-H\"{o}lder map. We start with the following Lemma.

\begin{lem} For any $x \in X\setminus M$ and any $\eps>0$, $f$
admits a $(1+\eps,\alpha)$-H\"{o}lder extension to $M \cup\{x\}$.
\end{lem}

\begin{proof} If $M$ is compact and $\delta >0$, we pick a $\delta$-net $\{x_1,..,x_n\}$ of
$M$ and a finite dimensional subspace $F$ of $Y$, containing
$f(x_1),..,f(x_n)$ such that $F$ is $(1+\delta)$-isomorphic to
some $\ell_\infty^m$. Then, there is $y\in F$ such that for all
$1\leq i\leq n$, $\Vert f(x_i)-y\Vert \leq (1+\delta)\Vert
x_i-x\Vert^\alpha$. If $\delta$ was chosen small enough, then for
any $z \in M$, $\Vert f(z)-y\Vert\leq (1+\eps)\Vert
z-x\Vert^\alpha.$

\noindent For a general $M$ and a fixed $x\in X\setminus M$, we
apply the compact case to the restriction of $f$ to $M\cap KB_X$,
for $K$ big enough and where $B_X$ denotes the closed unit ball of
$X$.

\end{proof}

We now finish the proof of Proposition \ref{almost}. Let
$(x_n)_{n\geq 1}$ be a dense sequence in $X\setminus M$. for a
given $\eps>0$, we pick $(\eps_n)_{n\geq 1}$ in $(0,1)$ so that
$\prod_{n\geq 1} (1+\eps_n)<1+\eps$. It follows from the above
Lemma and an easy induction that $f$ can be extended to a
$(1+\eps,\alpha)$-H\"{o}lder function on $M\cup \{x_n,\ n\geq 1\}$,
which in turn can be extended by density to $X$.

\end{proof}

\begin{rem} For $Y=C(K)$, there is a more concrete
argument, which even allows to extend $f$ isometrically when $M$
is compact. We use the Inf-convolution formula and define:

$$\forall t \in K\ f(x)(t)=\inf_{y\in M}[f(y)(t)+\Vert
x-y\Vert^\alpha].$$

\noindent Clearly, $\Vert f(x)-f(y)\Vert_\infty\leq \Vert
x-y\Vert^\alpha$. Since $f(M)$ is compact in $C(K)$, $f(x)$ is the
infimum of an equicontinuous family of functions and therefore is
continuous on $K$.
\end{rem}

\bigskip Let us now concentrate on the isometric extension
problem. We will need the following characterization.

\begin{lem}\lb{char} Let $(X,d)$ be a metric space, $M$ a
subset of $X$, $f:M\to c$ a contraction and $x \in X\setminus M$.
Then, the following statements are equivalent:
\begin{enumerate}
\item $f$ can be extended to a contraction $g:M\cup\{x\}\to c$.

\item $\forall \e>0 \quad \exists N\in{\bbN} \quad \forall n,m
> N \quad \forall y,z  \in M$
$$|f(y)(n)-f (z)(m)| \leq d(y,x)+d(z,x) + \e.$$
\end{enumerate}
\end{lem}

\begin{proof} Suppose that (1) holds. Then $(g(x)(n))_n$ is a
Cauchy sequence. Thus
$$\forall \eps>0\ \exists N\ \forall n,m>N\ \ \vert
g(x)(n)-g(x)(m)\vert<\eps.$$ Since $g$ is a contractive extension
of $f$, we have that for all $n,m>N$ and all $y,z \in M$

\begin{equation*}
\begin{split}
\vert f(y)(n)-f(z)(m)\vert &\le \vert f(y)(n)-g(x)(n)\vert+ \vert
g(x)(n)-g(x)(m)\vert+\vert g(x)(m)-f(z)(m)\vert \\ &\le
d(y,x)+d(z,x)+\eps.
\end{split}
\end{equation*}

Suppose now that (2) holds. Define
\begin{equation*}
s(j) = \sup_{m \geq j} \sup_{z \in M} (f(z)(m)-d(z,x)).
\end{equation*}

\noindent Let us fix $z_0 \in M$. Then, it is easily seen that
$$\forall j\in \Ndb,\ \ |s(j)|\leq \|f(z_0)\|_\infty+ d(x,z_0).$$
On the other hand $\{s(j)\}_{j \in \Ndb}$ is a decreasing sequence
and therefore converges. We will denote by $s(\infty)$ its limit.

\noindent In order to define $(g(x)(n))_{n\geq 1}$, we pick a
sequence $(N_k)_{k\geq 1}$ of integers such that

(i) (2) holds with $\e = 2^{-k}$ and $N=N_k$;

(ii) $\forall j > N_k \ \  s(j) \leq s(\infty)+2^{-k}$;

(iii) $\forall k\in \Ndb \ \ N_{k+1}> N_{k}$.

\noindent Then we define $g(x)$ as follows:

(1) for $n \leq N_1$, let $g(x)(n)$ be any element of
$$\bigcap_{y
\in M} [f(y)(n)-d(x,y),f(y)(n)+d(x,y)]=[\sup_{y\in
M}(f(y)(n)-d(x,y)),\inf_{y\in M}(f(y)(n)+d(x,y))].$$

(2) for $N_k < n\leq N_{k+1}$ we define
$$g(x)(n)=\max\{\sup_{y\in M}(f(y)(n)-d(x,y)),s(N_k)-2^{-k}\}.$$

\noindent It follows from (i) that
$$\forall n > N_k\ \forall y \in M\ \ s(N_k)-2^{-k}\leq
f(y)(n)+d(x,y).$$

\noindent So
$$\forall n\in \Ndb \ \ \sup_{y\in
M}(f(y)(n)-d(x,y))\leq g(x)(n)\leq \inf_{y\in
M}(f(y)(n)+d(x,y)).$$ Thus $g(x)\in \ell_\infty$ and for all $y$
in $M$, $\Vert g(x)-f(y)\Vert_\infty\leq d(x,y).$

\noindent Finally, note that
$$\forall n > N_k\ \ \sup_{y\in
M}(f(y)(n)-d(x,y)) \leq s(N_k).$$ Thus
$$\forall n \in (N_k,N_{k+1}]\ \ s(N_k)-2^{-k}\leq g(x)(n)\leq
s(N_k).$$ It is now clear that $(g(x)(n))_{n \geq 1}$ converges to
$s(\infty)$ and therefore belongs to $c$.
\end{proof}

As a first application we have

\begin{theorem} For any finite dimensional normed vector space $X$
$$(0,1)\subset {\cal A}(X,c).$$
\end{theorem}

\begin{proof} Let $0<\alpha<1$, $M$ a closed subset of $X$ such
that $0\notin M$ and $f:M\to c$ be a $(1,\alpha)$-H\"{o}lder function.
It is enough to show that $f$ admits a $(1,\alpha)$-H\"{o}lder
extension to $M\cup\{0\}$.

\noindent We fix $\eps>0$ and pick $x_0\in M$. Since $\alpha<1$,
$$\lim_{\Vert x\Vert \to \infty}[(\Vert x\Vert +\Vert
x_0\Vert)^\alpha- \Vert x\Vert^\alpha]=0.$$ So, there is $K>0$
such that $\Vert x-x_0\Vert^\alpha \leq \Vert x\Vert^\alpha
+\eps/3$ for all $x$ so that $\Vert x\Vert>K$. Let us also choose
$K$ such that $\Vert x_0\Vert \leq K$. Since $M_K=M\cap KB_X$ is
compact,

$$\exists N\in \Ndb \ \forall n,m>N \ \forall x \in M_K\ \
\vert f(x)(n)-f(x)(m)\vert< \frac{\eps}{3}.$$

\noindent Let now $x$ and $y$ in $M$.

\noindent If $x\in M_K$, then for all $n,m>N$:
$$\vert f(x)(n)-f(y)(m)\vert \leq \frac{\eps}{3} +\Vert x-y\Vert^\alpha
\leq \Vert x\Vert^\alpha +\Vert y\Vert^\alpha+\frac{\eps}{3}.$$ If
$x$ and $y$ belong to $M\setminus M_K$, then for all $n,m>N$:
$$\vert f(x)(n)-f(y)(m)\vert \leq \Vert x-x_0\Vert^\alpha +\Vert
y-x_0\Vert^\alpha +\frac{\eps}{3} \leq \Vert x\Vert^\alpha +\Vert
y\Vert^\alpha+\eps.$$ Then the conclusion follows directly from
Lemma \ref{char}.
\end{proof}

We will now see that the possibility of extending isometrically
all Lipschitz maps from a finite dimensional space into $c$ may
depend on the geometry of the space $X$. As a positive result, we
have for instance

\begin{theorem}\lb{linfty} For any $n \in \Ndb$
$${\cal A}(\ell_\infty^n,c)=(0,1].$$
\end{theorem}

\begin{proof} For $j\in\{1,\dots,n\}$,
$\delta\in \{-1,1\}$, we denote by $F_{j,\de}$ the following
$(n-1)-$face of the unit ball of $\ell_\infty^n$:
$$F_{j,\de}=\{x=(x_1,\dots,x_n) : \|x\|=1, x_j=\de\}.$$
Let $C_{j,\de}$ denote the cone supported by $F_{j,\de}$:
$$C_{j,\delta}=\{x\in\ell_\infty^n : x_j=\de\|x\|\}.$$
For $j,k\in \{1,..,n\}, j\ne k,$ and $\de,\eta\in\{-1,1\}$  we
denote by $F_{j,\de,k,\eta}$ the $(n-2)-$face of $F_{j,\de}$:
$$F_{j,\de,k,\eta}=F_{j,\de} \cap F_{k,\eta},$$
and by $C_{j,\de,k,\eta}$ the corresponding cone:
$$C_{j,\de,k,\eta}=C_{j,\de} \cap C_{k,\eta}.$$ We also define a family of
projections $P_{j,\delta,k,\eta}:C_{j,\de}\lra
C_{j,\delta,k,\eta}$ by
\begin{equation*}
P_{j,\delta,k,\eta}(x)=y,\ \text{where}\ \ \begin{cases} y_k=\eta|x_j| \\
y_i=x_i, \ \ \text{if $i\ne k$.}\end{cases}
\end{equation*}
Note that for every $x\in C_{j,\de}$, $\eta|x_j|=\eta\de x_j$, so
$P_{j,\delta,k,\eta}$ is linear on $C_{j,\de}$ and
\begin{equation}\lb{equal}
\forall x \in C_{j,\de}\ \ \|P_{j,\delta,k,\eta}(x)\|=\|x\|.
\end{equation}

\noindent Further, since for all $x\in C_{j,\de}$, $|x_j|\ge
|x_k|$ we get
\begin{equation}\lb{sign}
\sgn((P_{j,\delta,k,\eta}(x))_k-x_k)=\eta.
\end{equation}

\noindent We also introduce the projection $Q_k:\Rdb^n\to
\Rdb^{n-1}$ defined by
$$Q_k(x_1,..,x_n)=(x_1,..,x_{k-1},x_{k+1},..,x_n).$$

The following Lemma will provide us with a convenient finite
covering of the space $\ell_\infty^n$.

\begin{lem}\lb{partition}
For any $M\subset X=\ell_\infty^n$, any $\eps>0$ and any
$j\in\{1,\dots,n\}$, $\delta\in \{-1,1\}$, such that $C_{j,\delta}
\cap M \neq \emptyset$, there exist $A_1,..,A_\mu$ subsets of $X$
such that
$$(C_{j,\delta} \cap M)
\subset \bigcup_{i=1}^\mu A_i$$ and $\forall i \in\{1,..,\mu\}$
 $\exists x^i\in A_i\cap M$ satisfying
$$\forall x \in A_i\cap M \ \
\|x\|\geq \|x^i\|-\eps\ {\rm and}\  \Vert x-x^i\Vert \leq \Vert
x\Vert -\Vert x^i\Vert+\eps.$$
\end{lem}

\begin{proof}[Proof of Lemma~\ref{partition}] We will give a proof by induction on the dimension
of $\ell_\infty^n$. If $n=1$, the statement is clear, so let us
now assume that it is satisfied for $n-1$, where $n\geq 2.$

\noindent Let $M$, $\eps, j$ and $\delta$ be as in the statement
of Lemma \ref{partition}. We pick an element $x^{j,\delta} \in
C_{j,\delta}\cap M$ and we denote
$$B_{j,\delta}=x^{j,\delta}+C_{j,\delta}.$$
Note that
\begin{equation}\label{ineq1}
\forall x\in B_{j,\delta}\ \ \Vert x-x^{j,\delta}\Vert=\Vert
x\Vert -\Vert x^{j,\delta} \Vert.
\end{equation}

\noindent Denote
$d_{j,\delta,k,\eta}=|(P_{j,\delta,k,\eta}(x^{j,\delta}))_k -
(x^{j,\delta})_k|$. Let $x\in C_{j,\delta}$ such that for any
$k\in\{1,..,n\}\setminus\{j\},$ and any $\eta\in\{-1,1\}$,
\begin{equation}\label{big}
|(P_{j,\delta,k,\eta}(x))_k - x_k|\ge d_{j,\delta,k,\eta}
\end{equation}
Then, we claim that $x\in B_{j,\delta}$.

\noindent Indeed, by \eqref{sign}
$$|(P_{j,\delta,k,\eta}(x))_k - x_k|=\eta(P_{j,\delta,k,\eta}(x))_k - \eta x_k$$
and
$$|(P_{j,\delta,k,\eta}(x^{j,\delta}))_k -  (x^{j,\delta})_k|=
\eta(P_{j,\delta,k,\eta}(x^{j,\delta}))_k -
\eta(x^{j,\delta})_k.$$ Thus \eqref{big} implies that
$$\eta(P_{j,\delta,k,\eta}(x))_k-\eta(P_{j,\delta,k,\eta}(x^{j,\delta}))_k \ge \eta x_k- \eta(x^{j,\delta})_k,$$
and hence
$$\eta\eta\de x_j - \eta\eta\de (x^{j,\de})_j= \de(x_j-(x^{j,\de})_j)\ge \eta x_k- \eta(x^{j,\delta})_k.$$
Since this holds for all $\eta\in\{-1,1\}$, we get that for all
$k\in\{1,..,n\}\setminus\{j\}$,
$$\de(x-x^{j,\de})_j\ge |(x- x^{j,\delta})_k|.$$
Thus $x-x^{j,\de}\in C_{j,\delta}$ and $x\in B_{j,\delta}$.

\noindent Combining \eqref{sign} and \eqref{big}, we conclude that
$$C_{j,\delta}\setminus B_{j,\delta}\subset
\bigcup_{\substack{k\in\{1,..,n\}\setminus\{j\},\\
\eta\in\{-1,1\}}} B_{j,\delta,k,\eta}$$ where
$$B_{j,\delta,k,\eta}=\{x\in C_{j,\delta} : ((P_{j,\delta,k,\eta}(x))_k - x_k)\in \eta[0,d_{j,\delta,k,\eta})\}.$$

\noindent Now, for each $k \in\{1,..,n\}\setminus\{j\},$ and
$\eta\in\{-1,1\}$, we choose $N_{k,\eta}\in \Ndb$ such that
${\displaystyle \frac{d_{j,\delta,k,\eta}}{N_{k,\eta}}}<
\displaystyle{\frac{\eps}{3}}$. Then we set

\noindent $\forall k \in\{1,..,n\}\setminus\{j\}\ \forall
\eta\in\{-1,1\}\ \forall \nu \in \{1,\dots,N_{k,\eta}\}$:
$$I_{j,\delta,k,\eta}^\nu= \big[\frac{(\nu-1)d_{j,\delta,k,\eta}}
{N_{k,\eta}},\frac{\nu d_{j,\delta,k,\eta}}{N_{k,\eta}}\big)\ {\rm
and}\ B_{j,\delta,k,\eta}^\nu=\{x\in C_{j,\delta} :
((P_{j,\delta,k,\eta}(x))_k - x_k)\in \eta I_{j,\delta,k,\eta}^\nu
\}.$$ So we have
\begin{equation}\label{ineq2}
C_{j,\delta} \setminus
B_{j,\delta}=\bigcup_{\substack{k\in\{1,..,n\}\setminus\{j\},\\
\eta\in\{-1,1\}}}\bigcup_{\nu=1}^{N_{k,\eta}}
B_{j,\delta,k,\eta}^\nu.
\end{equation}

\noindent We now fix $k\in\{1,..,n\}\setminus\{j\},
\eta\in\{-1,1\}$ and $\nu\leq N_{k,\eta}$ such that
$B_{j,\delta,k,\eta}^\nu\cap M \neq \emptyset$ and denote for
simplicity:
$$B=B_{j,\delta,k,\eta}^\nu,\ I=\eta I_{j,\delta,k,\eta}^\nu ,
\ {\widetilde P}=P_{j,\de,k,\eta},\ P=Q_k{\widetilde P}, \
M'=P(M\cap B)$$
$$ {\rm and}\ C=P(C_{j,\delta})=
Q_kC_{j,\delta,k,\eta}=\{x\in\ell_\infty^{n-1} :
x_{\phi(j)}=\de\|x\|\},$$ where $\phi(j)=j$ if $k>j$ and
$\phi(j)=j-1$ if $k<j$.

\noindent Since $M'$ is a non empty subset of $C$, our induction
hypothesis yields the existence of $A'_1,..,A'_L \subset C$ so
that $M'\subset \bigcup_{l\leq L}A'_l$ and $\forall l \in
\{1,..,L\}$  $\exists y^l\in A'_l\cap M'$ satisfying
$$\forall y \in
A'_l\cap M'\ \ \|y\|\geq \|y^l\|-\frac\eps3\ {\rm and}\  \Vert
y-y^l\Vert \leq \Vert y\Vert -\Vert y^l\Vert+\frac\eps3.$$ Now let
$A_l=\{x=(x_i)_{i=1}^n\in C_{j,\delta},\ P(x)\in A_l', x_k \in
\de\eta x_j-I \}.$ We have that
$$B\cap M \subset \bigcup_{l\leq L}A_l.$$
Then, for any $l\leq L$, we pick $x^l\in A_l\cap M$ such that
$P(x^l)=y^l$. Note that
$$\forall x \in A_l\cap M\ \ \|x\|=|x_j|=\|Px\|\geq
\|y^l\|-\frac\eps3=|x_j^l|-\frac\eps3=\|x^l\|-\frac\eps3.$$
Therefore
$$\forall x \in A_l\cap M\ \ |x_j-x^l_j|\leq |x_j|-|x_j^l|+\frac{2\eps}3.$$
Now,
$$\Vert x-x^l\Vert
=\max\{\|P(x)-P(x^l)\|,|x_k-x_k^l|\}.$$ We have
$$\|P(x)-P(x^l)\|\leq
\|P(x)\|-\|P(x^l)\|+\frac\eps3=\|x\|-\|x^l\|+\frac\eps3.$$ Since
the diameter of $I$ is less than $\displaystyle{\frac\eps3}$, we
get on the other hand that
\begin{equation}
\begin{split}
|x_k-x_k^l|&= |(x_k-\eta\de x_j)-(x_k^l-\eta\de x_j^l)+\eta\de
x_j-\eta\de x_j^l| \\
&\leq \frac\eps3 +|x_j-x_j^l|\leq \eps
+|x_j|-|x_j^l|=\eps+\|x\|-\|x^l\|.
\end{split}
\end{equation}
So the conclusion of the lemma follows from \eqref{ineq1} and
\eqref{ineq2}.
\end{proof}

We now proceed with the proof of Theorem \ref{linfty}. As usual,
we consider a contraction $f:M\to c$, where $M$ is a closed subset
of $\ell_\infty^n$ with $0\notin M$. We will only show, as we may,
that $f$ can be contractively extended to $M\cup \{0\}$.

\noindent Let $\eps>0$. It follows from Lemma \ref{partition} that
there exist $A_1,..,A_\mu$ subsets of $X$ such that $M \subset
\bigcup_{i=1}^\mu A_i$ and
$$\forall 1\leq i\leq \mu\ \exists x^i\in A_i\cap M \ {\rm such\
that}\ \forall x \in A_i\cap M\  \ \Vert x-x^i\Vert \leq \Vert
x\Vert -\Vert x^i\Vert+\frac{\eps}{2}.$$ There also exists $N\in
\Ndb$ such that
$$\forall n,m>N\ \forall i\in\{1,..,\mu\}\ \ \vert
f(x^i)(n)-f(x^i)(m)\vert < \frac{\eps}{2}.$$ Let now $x$ and $y$
in $M$. Then we pick $i$ such that $x\in A_i$. Thus, for all
$n,m>N$
\begin{equation*}
\begin{split}
\vert f(x)(n)-f(y)(m)\vert &\le \vert f(x)(n)-f(x^i)(n)\vert+
\vert f(x^i)(n)-f(x^i)(m)\vert+\vert f(x^i)(m)-f(y)(m)\vert \\
&\le \Vert x-x^i\Vert+\frac{\eps}{2}+\|x^i-y\|\le \|x\|-\Vert x^i\Vert+\frac{\e}{2}+ \frac{\e}{2}+\Vert y\Vert +\|x^i\|\\
&\leq \Vert x\Vert+\Vert y\Vert+\eps.
\end{split}
\end{equation*}
Then we can apply Lemma \ref{char} to conclude our proof.

\end{proof}

\begin{cor}
Let $X$ be any finite dimensional Banach space whose unit ball is
a polytope. Then
$${\cal A}(X,c)=(0,1].$$
\end{cor}

\begin{proof} If $B_X$ is a polytope, we can find $f_1,...,f_n$ in the unit sphere of the
dual space of $X$ such that $$B_X=\bigcap_{i=1}^n \{x\in X,\
|f_i(x)|\leq 1\}.$$ Then  the map $T:X \to \ell_\infty^n$ defined
by $Tx=(f(x_i))_{i=1}^n$ is clearly a linear isometry and the
result follows immediately from Theorem~\ref{linfty}.
\end{proof}

We will finish this section with a counterexample in dimension 4.
We denote by $\ell_2^2\oplus_1\ell_2^2$ the space $\Rdb^4$
equipped with the norm:
$$\forall (s,t,u,v) \in \Rdb^4,\ \|(s,t,u,v)\|=(s^2+t^2)^{1/2}+(u^2+v^2)^{1/2}.$$
Then we have

\begin{theorem}\lb{counterex}  $${\cal
A}(\ell_2^2\oplus_1\ell_2^2,c)=(0,1).$$
\end{theorem}

\begin{proof} First we pick $K>1$ such that
\begin{equation}\label{ineq4}
\frac{1}{2}(\frac{K^2-1}{K^2})(\frac{K}{K+1})^3 > \frac{3}{8}.
\end{equation}

\noindent For $n \in \Ndb$, we define $x_n=(K^{2n},K^n,0,0)$ and
$y_n=(0,0,K^{2n},K^n)$. Note that
\begin{equation}\label{ineq5}
\forall n \in \Ndb,\ \ \Vert x_n\Vert \leq K^{2n}+\frac{1}{2}\ \
{\rm and}\ \ \Vert y_n\Vert \leq K^{2n}+{1 \over 2}.
\end{equation}
On the other hand,
$$\lim_{n \to \infty} (\Vert x_n\Vert - K^{2n})=
\lim_{n \to \infty} (\Vert y_n\Vert - K^{2n})=\frac{1}{2}.$$ So
\begin{equation}\label{ineq6}
\exists n_0\in \Ndb\ \forall n,m\geq n_0\ \ \Vert x_n-y_m\Vert
\geq K^{2n}+K^{2m}+\frac{7}{8}.
\end{equation}
Now, for all $n>m$, $K^n+K^m \leq K^n(\disp{\frac{K+1}{K}})$ and
$K^{2n}-K^{2m}\geq K^{2n}(\disp{\frac{K^2-1}{K^2}})$.

\noindent Since $$\Vert x_n-x_m\Vert
=(K^{2n}-K^{2m})[1+\frac{1}{(K^n+K^m)^2}]^{1/2},$$ we have
$$\Vert x_n-x_m\Vert\geq (K^{2n}-K^{2m})[1+(\frac{K}{K+1})^2
\frac{1}{K^{2n}}]^{1/2}.$$ Therefore, there exists $n_1\geq n_0$
such that for all $n>m\geq n_1$:
$$\Vert x_n-x_m\Vert \geq
(K^{2n}-K^{2m})[1+\frac{1}{2}(\frac{K}{K+1})^3 \frac{1}{K^{2n}}]
\geq
K^{2n}-K^{2m}+\frac{1}{2}(\frac{K^2-1}{K^2})(\frac{K}{K+1})^3.$$
Then, it follows from \eqref{ineq4} that for all $n>m\geq n_1$:
\begin{equation}\label{ineq7}
\Vert x_n-x_m\Vert \geq K^{2n}-K^{2m}+\frac{3}{8} \ \ {\rm and} \
\ \Vert y_n-y_m\Vert \geq K^{2n}-K^{2m}+\frac{3}{8}.
\end{equation}

Let us denote $M=\{x_n,\ n\geq n_1\}\cup\{y_n,\ n\geq n_1\}$. We
will now construct $u_n=f(x_n)$ and $v_n=f(y_n)$ in $c$ so that
$f:M\to c$ is 1-Lipschitz. So let $n\geq n_1$.

\noindent For $k$ odd and $k\leq n$, set
$u_n(k)=K^{2n}+\disp{\frac{5}{8}}$ and
$u_n(k)=K^{2n}+\disp{\frac{1}{4}}$ otherwise.

\noindent For $k$ even and $k\leq n$, set
$v_n(k)=-(K^{2n}+\disp{\frac{5}{8}})$ and
$v_n(k)=-(K^{2n}+\disp{\frac{1}{4}})$ otherwise.

We now check that $f$ is 1-Lipschitz.

\noindent For all $n>m\geq n_1$, $\Vert u_n-u_m\Vert_\infty\leq
K^{2n}+\disp{\frac{5}{8}}
-(K^{2m}+\disp{\frac{1}{4}})=K^{2n}-K^{2m}+\disp{\frac{3}{8}}$.

\noindent Therefore, by \eqref{ineq7}, $\Vert u_n-u_m\Vert_\infty
\leq \Vert x_n-x_m\Vert$.

\noindent We have, as well that $\Vert v_n-v_m\Vert_\infty \leq
\Vert y_n-y_m\Vert$.

\noindent We also have that for all $n,m\geq n_1$, $\Vert
u_n-v_m\Vert_\infty=K^{2n}+K^{2m}+\disp{\frac{7}{8}}$.

\noindent Thus, \eqref{ineq6} implies that $\Vert
u_n-u_m\Vert_\infty \leq \Vert x_n-y_m\Vert$.

\noindent We have shown that $f$ is 1-Lipschitz.

Assume now that $f$ can be extented at 0 into a 1-Lipschitz
function $g$ and let $g(0)=w=(w(k))_{k\geq 1} \in c$. Then it
follows from \eqref{ineq5} that for all odd values of $k$,
$w(k)\geq \disp{\frac{1}{8}}$ and for all even values of $k$ $w(k)
\leq -\disp{\frac{1}{8}}$. This contradicts the fact that $w \in
c$.

\end{proof}

\begin{rem}  As we already mentioned in the
introduction, this seems to be the first example of Banach spaces
$X$ and $Y$ such that ${\cal A}(X,Y)\neq \widetilde{\cal A}(X,Y)$
and also such that ${\cal A}(X,Y)$ is not closed in $(0,1]$.
\end{rem}

\section{Maps into $C(K)$ spaces}

In this last section we show that if $K$ is an infinite compact
metric space, then the study of the isometric extension for
Lipschitz maps with values in $C(K)$ reduces to the results of the
previous section. More precisely, we   prove the following.

\begin{theorem}\lb{C(K)} Let $(X,d)$ be a metric space and $(K,\ro)$ be an
infinite compact metric space. Then
$${\cal A}(X,C(K))={\cal A}(X,c).$$
\end{theorem}

The main step of the proof will be to establish the following
generalization of Lemma \ref{char}.

\begin{prop}\lb{Char} Let $M$ be a
subset of $X$, $f:M\to C(K)$ a contraction and $x \in X\setminus
M$. We denote by $D$ the diameter of $K$ for the distance $\ro$.
Then, the following statements are equivalent:
\begin{enumerate}
\item $f$ can be extended to a contraction $g:M\cup\{x\}\to C(K)$.

\item $\forall \e>0 \quad \exists \de>0$ such that $\forall t,s\in K$
with $\ro(t,s)<\de  \quad \forall y,z  \in M$
$$|f(y)(t)-f (z)(s)| \leq d(y,x)+d(z,x) + \e.$$

\item $\exists \f:[0,D]\lra [0,+\infty)$ such that $\f$ is continuous, $\f(0)=0$ and
$$\forall t,s\in K \ \forall y,z  \in M \ \ |f(y)(t)-f (z)(s)|
\leq d(y,x)+d(z,x) + \f(\ro(t,s)).$$

\end{enumerate}
\end{prop}

\begin{proof} Suppose that (1) holds. Then (2) follows from the triangle
inequality and the fact that $g(x)$ is uniformly continuous on
$K$.

Assume now that (2) holds. Let us define, for $\lambda \in (0,D]$:
$$\xi(\lambda)=\sup_{y,z\in M}\sup_{\ro(t,s)\leq
\lambda}(|f(y)(t)-f(z)(s)|-d(x,y)-d(x,z)).$$ The function $\xi$ is
clearly non decreasing and bounded below by $-2\,dist\,(x,M)$. So
we can set
$$\xi(0)= \lim_{\lambda\searrow 0}\xi (\lambda).$$
We have that
$$\forall t,s \in K\ \forall y,z \in M\ \ |f(y)(t)-f (z)(s)|
\leq d(y,x)+d(z,x) + \xi(\ro(t,s)).$$ It follows from (2) that
$\xi(0)\leq 0$. So, if we set $\psi=\xi -\xi(0)$, we get that
$\psi$ is non decreasing, $\psi(0)=0$ and $\psi$ is continuous at
0. Since $\psi\ge\xi$, we still have
$$\forall t,s \in K\ \forall y,z \in M\ \ |f(y)(t)-f (z)(s)|
\leq d(y,x)+d(z,x) + \psi(\ro(t,s)).$$
 We now define the function $\f$ in the following way:
$\f (0)=0$ and for $n\in \Ndb$, $\f(\frac D{n+1})=\psi(\frac
D{n})$. We also ask $\f$ to be constant equal to $\psi(D)$ on
$[\frac D2,D]$, and affine on each $[\frac D{n+2},\frac D{n+1}]$
for $n\in \Ndb$. It is now clear that $\f$ is non decreasing,
continuous on $[0,D]$ and that $\psi\leq \f$ on $[0,D]$. So we
have
$$\forall t,s \in K\ \forall y,z \in M\ \ |f(y)(t)-f (z)(s)|
\leq d(y,x)+d(z,x) + \f(\ro(t,s)).$$ This proves that (2) implies
(3).

Suppose now that (3) holds and define, for $t \in K$,
\begin{equation*}
g(x)(t) = \sup_{s\in K} \sup_{z \in M}
(f(z)(s)-d(z,x)-\f(\ro(t,s))).
\end{equation*}

\noindent Fix $y_0\in M$. Then, for all $z\in M$ and for all $s\in
K$,
$$f (z)(s)- d(z,x) - \f(\ro(t,s))\leq \|f(y_0)\|_{C(K)}+d(z,y_0)-d(z,x)
\leq \|f(y_0)\|_{C(K)}+d(x,y_0).$$ So $g(x)(t)$ is well defined.
Further, it follows from the uniform continuity of $\f$ on $[0,D]$
that $g(x)$ is continuous on $K$.

\noindent Since $\f(0)=0$, we have, by definition of $g(x)$, that
for all $y\in M$ and all $t\in K$
\begin{equation}\lb{ck1}
f(y)(t)-g(x)(t)\le d(x,y).
\end{equation}
By (3), we get that for all $y,z\in M$ and for all $t,s\in K$
$$ |f(z)(s)-f (y)(t)| \leq d(y,x)+d(z,x) + \f(\ro(t,s)),$$
so
$$f (z)(s)- d(z,x) - \f(\ro(t,s))\leq f(y)(t)+d(y,x),$$
and by taking the supremum over $z$ and $s$ we obtain
\begin{equation}\lb{ck2}
g(x)(t)-f(y)(t)\le d(x,y).
\end{equation}
Combining \eqref{ck1} and \eqref{ck2}, we get that for all $y\in
M$ $\|g(x)-f(y)\|_{C(K)}\le d(x,y)$. Thus (3) implies (1) and this
ends the proof of Proposition \ref{Char}.

\end{proof}

\begin{proof}[Proof of Theorem \ref{C(K)}] Since $K$ is an infinite
compact metric space, it contains a closed subset $F$ which is
homeomorphic to the one point compactification of $\Ndb$. Then,
$C(F)$ is clearly isometric to $c$. On the other hand, by the
linear version of Tietze extension theorem due to K. Borsuk
\cite{Bo}, there is a linear isometry $T:C(F) \to C(K)$ such that
for any $f$ in $C(F)$, $Tf$ is an extension of $f$ to $K$. Let now
$R$ be the restriction operator from $C(K)$ onto $C(F)$. Then
$P=TR$ is a projection of norm 1 from $C(K)$ onto an isometric
copy of $c$. Therefore, it is clear that for any metric space $X$,
${\cal A}(X,C(K))\subset {\cal A}(X,c)$.

For the other inclusion, it is enough to show that if $1\notin
{\cal A}(X,C(K))$, then $1 \notin {\cal A}(X,c)$. So let us assume
that $1\notin {\cal A}(X,C(K))$. Then there exist $M \subset X$, a
contraction $f: M\to C(K)$ and $x\in X\setminus M$ such that $f$
can not be contractively extended to $M \cup \{x\}$. Thus, by
Proposition~\ref{Char}, there exists $\e>0$ so that for all
$n\in\bbN$ there exist $t_n, s_n\in K$ with $\ro(t_n,s_n)<1/n$ and
$y_n, z_n\in M$ so that
\begin{equation}\lb{equal1}
|f(y_n)(t_n)-f (z_n)(s_n)| > d(y_n,x)+d(z_n,x) + \e.
\end{equation}

\noindent Since $K$ is compact, we may assume that the sequence
$(t_n)_{n\in\bbN}$ is convergent. Define now a sequence
$(w_n)_{n\in\bbN}$ in $K$ by setting, for $n\in \Ndb$,
$w_{2n-1}=t_n$ and $w_{2n}=s_n$. Then the sequence
$(w_n)_{n\in\bbN}$ is convergent. So we can define a 1-Lipschitz
map $h :M \to c$ by
$$\forall y \in M\ h(y)=(h(y)(n))_{n\in \Ndb}=(f(y)(w_n)))_{n\in
\Ndb}.$$ It now clearly follows from \eqref{equal1} and
Lemma~\ref{char} that $h$ does not have any extension to a
1-Lipschitz map from $M\cup\{x\}$ into $c$. Therefore $1\notin
{\cal A}(X,c)$.

\end{proof}

\end{document}